# Higher ρ Invariants
## Shmuel Weinberger[1]

to Paul Baum on his 60<sup>th</sup> birthday

The ρ-invariant is an invariant of odd dimensional manifolds with finite fundamental group, and lies in the representations modulo the regular representations ($\otimes \mathbb{Q}$). It is a fundamental invariant that occurs in classifying lens spaces [AB], their homotopy analogues [W1], and is intimately related to η-invariant for the signature operator [APS]. The goal of this note is to use some of the technology developed in studying the Novikov higher signature conjecture to define an analogous invariant for certain situations with infinite fundamental group.

For motivation, let us consider the problem of distinguishing various particular manifolds. For openers, consider the 3-dimensional lens spaces $L_5(1,1)$ and $L_5(2,3)$; the notation here means the lens spaces with fundamental group $\mathbb{Z}_5$ and with rotation numbers ($2\pi/5$, $2\pi/5$) in the first case and ($4\pi/5$, $6\pi/5$) in the second. These lens spaces are homotopy equivalent (the product of the rotation numbers determines this: 1×1 = 2×3 mod 5). There are essentially two classical proofs that these are not diffeomorphic. The most classical uses Reidemeister torsion and can be found in [Mi] or [Co]. The second proof is given in [AB] and reformulated analytically in [APS] and in a topological form in [W1] (with some more topologization in [Gi]).

Unfortunately, neither proof directly applies to the products $S^1 \times L_5(1,1)$ versus $S^1 \times L_5(2,3)$. One would not expect the torsion proof to work in any form because, although one can use torsion to distinguish h-cobordant manifolds from one another, these manifolds do become diffeomrophic after crossing with the circle!

However, one can deduce the result indirectly from the Atiyah-Bott proof. If these manifolds were diffeomorphic, by passing to the infinite cyclic cover, one could embed $L_5(1,1)$ in $\mathbb{R}^1 \times L_5(2,3)$, and thus obtain an h-cobordism between $L_5(1,1)$ and $L_5(2,3)$, which is, in fact, what [AB] show can't occur.

However, once we cross with a second circle, this trick does not suffice. Nonetheless, one can indirectly prove that these remain different using Shaneson's thesis [Sh]: aside from 2-torsion (of which there is none in this lens space) $S(M) \dashrightarrow S(M \times S^1)$ is always injective, when $S(M)$ denotes the manifolds simple homotopy equivalent to M, up to diffeomorphism. Our goal is then to find invariants that distinguish these manifolds.

Now, one would be tempted to try to invoke "higher" analogues of the [AB, APS] invariant analogous to what does with the ordinary signature. However, there is a key difference between signature and eta type invariants: the former is cobordism invariant while the latter is not: indeed the latter measures signatures of nulcobordisms.

Perhaps more strikingly: although the manifolds $S^1 \times L_5(1,1)$ and $S^1 \times L_5(2,3)$ are nondiffeomorphic, one can readily deduce from [CS] that after taking connected sums with sufficiently many copies of $S^2 \times S^2$ they become diffeomorphic! (Note that the diffeomrphism is <u>not</u> homotopic to the original homotopy equivalence $\#$ id -- it necessarily moves around the middle dimensional homology.)

This example shows that, unlike the classical situation of odd dimensional manifolds with finite fundamental group, one cannot expect to be able to define the invariant for all manifolds of a given dimension with a given $\pi_1$.

Definition: M is said to be <u>antisimple</u> if the chain complex of $M^n$ is chain equivalent to complex of projective modules $P_*$ with $P_i=0$ for i in the middle dimension.

This is a slight variant of a definition of Hausmann, who used a geometric form of this: namely, that M have a handlebody with no middle dimensional handles. These conditions are closely related, with the algebraic version given above being slightly weaker. It can be described purely cohomologically as in [W2] as that $H^i(M; \Lambda) = 0$ for all projective $\mathbb{Q}\pi$ modules, $\Lambda$, for $i = [n/2]$.

The key example to have in mind is $S^1 \times L_5(2,3) \# S^2 \times S^2$ which clearly is not antisimple. Lens spaces are antisimple if one uses $\mathbb{Q}\pi$ coefficients. (One needs some little trickery to make tori × lens spaces antisimple for high dimensional tori. I leave it to the reader to invent her own fix.)

Theorem: Suppose M is an antisimple manifold with fundamental group $\Gamma$. Suppose that the Novikov conjecture is true for $\Gamma$. Then one can associate to M an element of
$$H\rho(M) \in L_{n+1}(\Gamma) \otimes \mathbb{Q} / \oplus_{i \geq 0} H_{n+1-4i}(B\Gamma; \mathbb{Q}).$$
The invariant ranges over a lattice in this vector space as one varies among manifolds homotopy equivalent to M.

Remark: We will see that this invariant measures, in some way, the depth of the Novikov conjecture: that is, our invariants are secondary η-type invariants associated to the solutions of the Novikov conjecture.

<u>Remark</u>: A solution to the proper equivariant Novikov for $\Gamma$ (see [FRW, RsW, BC, BCH]) then gives a map from $L_{n+1}(\Gamma)$ to an equivariant K-group. The latter has a map to a sum of homology groups, parametrized by the conjugacy classes of elements of finite order in $\Gamma$ (see [BC]). For instance for a group of the form $\Gamma \times \mathbb{Z}_2$, where $\Gamma$ is torsion free, one gets two copies of $\oplus H_{n+1\pm 4i}(B\Gamma;\mathbb{Q})$. (If we'd repalced $\mathbb{Z}_2$ by $\mathbb{Z}_3$ we'd get $\oplus H_{n+1\pm 4i}(B\Gamma;\mathbb{Q}) \times \oplus H_{n+1\pm 2i}(B\Gamma;\mathbb{Q})$.) The homology piece that we have to mod out in the above formula maps into the part that corresponds to the identity element. In other words, for the nontrivial elliptic elements, one has no indeterminacy in the invariant, i.e. an invariant in $\oplus H_{n+1\pm 4i}(B\Gamma;\mathbb{Q})$ (or $\oplus H_{n+1\pm 2i}(B\Gamma;\mathbb{Q})$) and for the trivial element, one gets an invariant in $\oplus H_{n+4i+5}(B\Gamma;\mathbb{Q})$.

One only needs a solution to the usual Novikov conjecture to get the invariant lying in $\oplus H_{n+4i+5}(B\Gamma;\mathbb{Q})$, so for an extremely large clas of groups one obtains this type of invariant.

Examples of groups to which this theory can be applied include discrete subgroups of Lie groups, hyperbolic groups, etc. (see the proceedings of the Oberwolfach conference on the Novikov Conjecture for a reasonably recent summary.)

<u>Remark</u>: Some other settings where one can define Hρ were the content of an early version of this paper, referred to in [Lo], and will be briefly described below. A special case of those invariants was tacitly referred to in one of the corollaries in [We1] to describe bordism of homologically trivial group actions. The most novel aspect of the current version is the piece associated to the identity elliptic element. The need for this part arose in joint work with Alex Nabutovsky on the undecidability of homeomorphism for manifolds with <u>given fundamental group</u>, and will be discussed in a companion paper to this one.

<u>Proof of theorem</u>: We will show that in appropriate senses M bounds two different objects, which when glued together give the element of $L_{n+1}$. The homology term is the indeterminacy.

By hypothesis, $C_*(M)$ is chain equivalent to $P_*$ where $P_i = 0$ for $i = [n/2]$. Consider the pair of complexes $(C_*, P^{<[n/2]})$ where $P^{<[n/2]}$ is just $P_*$ truncated at $[n/2]$. One can check that this complex natually can be given the structure of an algebraic Poincare pair in the sense of [R]. This is our first boundary[2].

Our second boundary is more geometric, and follows a trick I used in [We2] to study ordinary $\eta$. Of course the mere fact that $\sigma^*(M) = 0$ in $L^n(\Gamma)$ is not enough for us to conclude that M bounds a manifold with fundamental group $\Gamma$, but when the Novikov conjecture holds for $\Gamma$, we can get quite close.

Recall that, following [S], a Witt space is an oriented polyhedral pseudomanifold whose intersection chain (with middle perversity) sheaf is a self dual sheaf[3]. Interestingly, the bordism of Witt spaces forms an extermely nice homology theory, whose rational calculation is simply:

$$\Omega_n^{Witt}(X) \otimes \mathbb{Q} \cong \oplus_{i \geq 0} H_{n-4i}(X; \mathbb{Q}).$$

---

[2]A similar point is made in [DR]; our motivation was the Deligne construction of the intersection chain sheaf.
[3]The definition in [S] is equivalent, but it does not really matter. We can use his definition just as well, and rely onthe fact that he proves that IC is self dual for (his) Witt spaces.

As observed in [CSW] one can define symmetric signatures for Witt spaces, and there is a commutative diagram

$$\begin{array}{ccc} \Omega_n(B\Gamma) & \longrightarrow & \Omega_n^{Wit}(B\Gamma) \\ \downarrow & & \downarrow \\ \oplus_{i \geq 0} H_{n-4i}(B\Gamma;\mathbb{Q}) & \longrightarrow & L^n(\Gamma) \otimes \mathbb{Q} \end{array}$$

In fact, the bottom left and upper right groups are naturally isomorphic, using the L-classes produced by the method of Thom-Milnor and the cobordism invariant signature of Witt spaces. (See [GM, S].) Thus, the Novikov conjecture, which asserts te bottom arrow to be injective implies that $\Omega_n^{Wit}(B\Gamma) \longrightarrow L^n(\Gamma) \otimes \mathbb{Q}$ is injective, and that therefore, in our situation, some multiple of M, mM bounds a Witt space, we shall deonte by W, such that $\pi_1(M) \longrightarrow \pi_1(W)$ is an isomorphism.

Now $H\rho(M) = m^{-1} [IC^m_*(W) \cup_m P^{<[n/2]}] \in L^{n+1}(\Gamma) \otimes \mathbb{Q}$ where "union" is taken in the sense of Ranicki's glueing of algebraic Poincare complexes [R]. Now the only question we have is the ambiguities of this construction.

Note that $(A \cup B) \cup -(A' \cup B') = (A \cup -A') \oplus (B \cup -B')$

The indeterminacy in the choice of P is only up to chain autohomotopiy equivalence. $P \cup -P'$ will be another closed antisimple Poincare space, and therefore has vanishing symmetric signature. On the other hand, the indeterminacy in the Witt space is a closed Witt space of dimension n+1, which can be identified with $\oplus_{i \geq 0} H_{n+1-4i}(B\Gamma;\mathbb{Q})$.

The realization of values is an exercise in surgery exact sequences, Wall realization and assembly maps. We leave it to the reader.
    QED

Remark: With more work, it is possible to get some integral information; we leave this for an interested reader.

Remark: One can easily deduce multiplicativity formulae for Hρ, at least when it makes sense. (Products are not ordinarily antisimple.)

Remark: It seems reasonable to expect an analytic interpretation of this invariant a la higher eta invariants which would apply these to situations where elliptic operators are known to be invertible. For instance, presumably one can use such things to distinguish components in spaces of positive scalar curvature metrics. In particular, one would like to know whether the different components of the moduli space of positive scalar curvature metrics on the sphere remain distinct after crossing with some other manifold. One might well conjecture that crossing with an aspherical manifold does not allow the connection of these components[4]. This particular result can often be verified directly (at least for components detected by a relative index theorem). [Lo] achieves the whole analytic program when the fundamental group has polynomial growth.

Final Remark: The remaining contexts where Hρ was defined in the old version of this paper were two fold:

---

[4]This includes the well-known conjecture that aspherical manifolds do not themselves possess metrics of positive scalar curvature.

(a) Homologically trivial group actions in the sense of [We1]. In that case the relevant fundamental group is $\Gamma \times \pi$, where $\pi$ is finite, and one has acyclicity over the ring $\mathfrak{I}\Gamma$, where $\mathfrak{I}$ is the augmentation ideal of $\mathbb{Q}\pi$. Then the invariant lies in $H\rho(M) \in L_{n+1}(\mathfrak{I}\Gamma)/\oplus_{i \geq 0} H_{n+1-4i}(B\Gamma;\mathbb{Q})$. The analogue of the vanishing result for antisimple spaces necessary here is a result from [We3] that gives vanishing of higher signatures for manifolds with free homologically trivial group actions, assuming the Novikov conjecture. As mentioned above, this version plays a role in the cobordism classification results in [We1]

(b) More tautologously, given a (rational) homotopy equivalence M' ---> M, so that one has a structure, then one can use the contravariant functoriality of surgery (see [R2, We5]) to define an element in $S(B\Gamma; \mathbb{Q}) \cong L_{n+1}(\mathbb{Q}\Gamma)/\oplus H_{n+1 \pm 4i}(B\Gamma;\mathbb{Q})$. Most of the invariant comes from the nontrivial elliptic elements in a (conjectural) cyclic homology description of L. The "new piece" of the invariant associated to the the identity elliptic element arises in this thery as follows (following on the discussion in the beginning of [We4]). A <u>solution</u> to the Novikov conjecture gives a map $S(M)$ ---> $H_{n+1}(B\Gamma, M; L)$, which refines the difference of the L-classes of M and M'. Since M is n-dimensional we can map this group to the quotient of $H_{n+1}(B\Gamma; L)$ by the image of the homology of the n skeleton of $B\Gamma$.

(c) In a future paper with Misha Farber we will describe a theory of higher signature for manifolds with boundary, even if the latter is not antisimple. A key feature of that situation is that one must have a failure of Novikov additivity in the nonsimply connected case.

Acknowledgements: I would like to thank John Lott and Alex Nabutovsky for useful conversations. I'd

especially like to thank Lott for suggesting that one look at signature operators on manifolds with boundary where the operator on the boundary is already invertible as a variant of the homologically trivial actions considered in the first version of this paper. This condition seems to be closely related to our antisimplicity hypothesis. I would also like to mention the forthcoming University of Maryland thesis of Navin Keswani which contains some interesting related material.